\documentclass[12pt,reqno]{amsart}

\usepackage{amsfonts,amssymb,amsmath}

\DeclareMathOperator{\bS}{\boldsymbol{\mathfrak{S}}}
\DeclareMathOperator{\bF}{\boldsymbol{\mathfrak{F}}}
\DeclareMathOperator{\bH}{\boldsymbol{\mathfrak{H}}}

\DeclareMathOperator{\size}{\mathit{d}}
\DeclareMathOperator{\lin}{\mathrm{lin}}
\DeclareMathOperator{\cone}{\mathrm{pos}}
\DeclareMathOperator{\conv}{\mathrm{conv}}

\newtheorem{theorem}{Theorem}

\newtheorem{corollary}[theorem]{Corollary}

\newtheorem{proposition}[theorem]{Proposition}

\newtheorem{example}[theorem]{Example}

\numberwithin{equation}{section} \numberwithin{theorem}{section}

\begin{document}

\title{Faces and Bases: Dehn-Sommerville Type Relations}

\date{}
\author{Andrey O. Matveev}
\address{Data-Center Co., RU-620034, Ekaterinburg,
P.O.~Box~5, Russian~Federation} \email{aomatveev@\{dc.ru,
hotmail.com\}}


\begin{abstract}
We review several linear algebraic aspects of the Dehn-Sommerville
relations and relate redundant analogues of the $f$- and
$h$-vectors describing the subsets of a simplex
$\mathbf{2}^{\{1,\ldots,m\}}$ that satisfy Dehn-Sommerville type
relations to integer points contained in some rational polytopes.
\end{abstract}

\maketitle

\section{Introduction and preliminaries}
\thispagestyle{empty}

A family $\Delta$ of subsets of a finite set $V$ is called an {\em
abstract simplicial complex\/} (or a {\em complex}) on the {\em
vertex\/} set $V$ if the inclusion $A\subseteq B\in\Delta$ implies
$A\in\Delta$, and if $\{v\}\in\Delta$, for any $v\in V$; see,
e.g.,~\cite{BB,B,BH,BP,Hibi,MS,St1,Z}. The sets from the family
$\Delta$ are called the {\em faces}. The {\em size\/} of a face
$F\in\Delta$ is its cardinality $|F|$. Let $\#$ denote the number
of sets in a family. If $\#\Delta>0$ then the {\em size\/}
$\size(\Delta)$ of $\Delta$ is defined by
$\size(\Delta):=\max_{F\in\Delta}|F|$.

The row vector
$\pmb{f}(\Delta):=\bigl(f_0(\Delta),f_1(\Delta),\ldots,
f_{\size(\Delta)-1}(\Delta)\bigr) \in\mathbb{N}^{\size(\Delta)}$,
where $f_i(\Delta):=\#\{F\in\Delta:\ |F|=i+1\}$, is called the
{\em $f$-vector\/} of $\Delta$;  these vectors are characterized
by the Sch\"{u}tzenberger-Kruskal-Katona theorem. The row {\em
$h$-vector}
$\pmb{h}(\Delta):=\bigl(h_0(\Delta),h_1(\Delta),\ldots,h_{\size(\Delta)}\bigr)
\in\mathbb{Z}^{\size(\Delta)+1}$ of $\Delta$ is defined by
\begin{equation*}
\sum_{i=0}^{\size(\Delta)}h_i(\Delta)\cdot\mathrm{y}^{\size(\Delta)-i}:=
\sum_{i=0}^{\size(\Delta)}f_{i-1}(\Delta)\cdot(\mathrm{y}-1)^{\size(\Delta)-i}\
.
\end{equation*}

For a positive integer $m$, let $[m]$ denote the set
$\{1,2,\ldots,m\}$. We denote by $\mathbf{2}^{[m]}$ the {\em
simplex\/} $\{F:\ F\subseteq[m]\}$ and relate to an arbitrary {\em
face system\/} $\Phi\subseteq\mathbf{2}^{[m]}$ ``long'' analogues
of the $f$- and $h$-vectors, namely, the row vectors
\begin{align*}
\pmb{f}(\Phi;m):&=\bigl(f_0(\Phi;m),f_1(\Phi;m),\ldots,f_m(\Phi;m)\bigr)
\in\mathbb{N}^{m+1}\ ,\\
\pmb{h}(\Phi;m):&=\bigl(h_0(\Phi;m),h_1(\Phi;m),\ldots,h_m(\Phi;m)\bigr)
\in\mathbb{Z}^{m+1}\ ,
\end{align*}
where $f_i(\Phi;m):=\#\{F\in\Phi:\ |F|=i\}$, for $0\leq i\leq m$,
and the vector $\pmb{h}(\Phi;m)$ is defined by
\begin{equation*}
\sum_{i=0}^m h_i(\Phi;m)\cdot\mathrm{y}^{m-i}:=\sum_{i=0}^m
f_i(\Phi;m)\cdot(\mathrm{y}-1)^{m-i}\ ;
\end{equation*}
see, e.g.,~\cite{McMSh,McMWa}, where similar combinatorial tools
appear. The maps $\Phi\mapsto\pmb{f}(\Phi;m)$ and
$\Phi\mapsto\pmb{h}(\Phi;m)$ from the Boolean lattice
$\mathcal{D}(m)$ of all face systems (ordered by inclusion) to
$\mathbb{Z}^{m+1}$ are {\em valuations\/} on $\mathcal{D}(m)$. We
consider the vectors $\pmb{f}(\Phi;m)$ and $\pmb{h}(\Phi;m)$ as
elements from the real Euclidean space $\mathbb{R}^{m+1}$ of row
vectors.

Given a complex $\Delta\subseteq\mathbf{2}^{[m]}$ with
$\#\Delta>0$, denote by $\mathbf{U}\bigl(\size(\Delta)\bigr)$ the
$(\size(\Delta)+1)\times(\size(\Delta)+1)$ {\em backward identity
matrix} whose rows and columns are indexed starting with zero and
whose $(i,j)$th entry is the Kronecker delta
$\delta_{i+j,\size(\Delta)}$. One says that the $h$-vector of a
complex $\Delta$ satisfies the {\em Dehn-Sommerville relations\/}
if $\pmb{h}(\Delta)$ is a left eigenvector of
$\mathbf{U}\bigl(\size(\Delta)\bigr)$ corresponding to the
eigenvalue $1$, that is, it holds
\begin{equation}
\label{eq:12} h_l(\Delta)=h_{\size(\Delta)-l}(\Delta)\ ,\ \ \
0\leq l\leq\size(\Delta)
\end{equation}
(see, e.g.,~\cite[\S{}VI.6]{BarvinokConv}, \cite[Chapter~5]{BR},
\cite[\S{}II.5]{BH}, \cite[\S\S{}1.2, 3.6, 8.6]{BP},
\cite[\S{}III.11]{Hibi}, \cite[\S\S{}II.3, II.6, III.6]{St1},
\cite[\S{}3.14]{St2}, \cite[\S{}8.3]{Z},) or
\begin{equation}
\label{eq:1}
h_l(\Delta;m)=(-1)^{m-\size(\Delta)}h_{m-l}(\Delta;m)\ ,\ \ \
0\leq l\leq m\ ,
\end{equation}
see, e.g.,~\cite[p.~171]{McMSh}.

In this paper we consider the set of all row vectors
$\mathfrak{f}\in\mathbb{N}^{m+1}$ such that for each of them there
exists a face system $\Phi\subset\mathbf{2}^{[m]}$ satisfying
Dehn-Sommerville type relations analogous to (\ref{eq:1}), and
$\pmb{f}(\Phi;m)=\mathfrak{f}$; we interpret the set of those
vectors $\mathfrak{f}$ as a set of integer points contained in two
rational polytopes.

\section{Notation}
\label{section:notation}

Throughout the paper, $m$ denotes an integer, $m\geq 2$.

The components of vectors as well as the rows and columns of
matrices are indexed starting with zero.

$\mathbf{I}(m)$ denotes the $(m+1)\times(m+1)$ {\em identity
matrix}.

$\mathbf{T}(m)$ is the {\em forward shift matrix} whose $(i,j)$th
entry is $\delta_{j-i,1}$.

For a vector $\pmb{w}:=(w_0,\ldots,w_m)$, $\pmb{w}^{\top}$ denotes
its transpose; $\|\pmb{w}\|^2:=\sum_{i=0}^m w_i^2$.

If $\boldsymbol{\mathfrak{B}}:=(\pmb{b}_0,\ldots,\pmb{b}_m)$ is a
basis of $\mathbb{R}^{m+1}$ then, given a vector
$\pmb{w}\in\mathbb{R}^{m+1}$, we denote by
$[\pmb{w}]_{\boldsymbol{\mathfrak{B}}}:=\bigl(
\kappa_0(\pmb{w},\boldsymbol{\mathfrak{B}}),\ldots,
\kappa_m(\pmb{w},\boldsymbol{\mathfrak{B}})\bigr)\in\mathbb{R}^{m+1}$
the $(m+1)$-tuple satisfying the equality
$\sum_{i=0}^m\kappa_i(\pmb{w},\boldsymbol{\mathfrak{B}})\cdot\pmb{b}_i=\pmb{w}$.
We also set $\kappa_i(\pmb{w}):=\kappa_i(\pmb{w},\bS_m)=w_i$,
$0\leq i\leq m$; $\bS_m$ stands for the standard basis of
$\mathbb{R}^{m+1}$.

$\pmb{0}:=(0,0,\ldots,0)$;
$\boldsymbol{\iota}(m):=(1,1,\ldots,1)$.

$\lin(\cdot)$, $\cone(\cdot)$ and $\conv(\cdot)$ stand for a
linear, conical, and convex hulls, respectively.

$\pmb{V}\oplus\pmb{W}$ denotes the direct sum of subspaces, and
$\boxplus$ denotes the Minkowski addition.

We use the notation $\hat{0}$ to denote the empty set.

If $A\subseteq C\in\mathbf{2}^{[m]}$ then $[A,C]$ denotes the {\em
Boolean interval\/} $\{B\in\mathbf{2}^{[m]}:\ A\subseteq
B\subseteq C\}$.

If $\Phi$ is a face system, $\#\Phi>0$, then its {\em size\/}
$\size(\Phi)$ is defined by $\size(\Phi):=\max_{F\in\Phi}|F|$.

\section{Dehn-Sommerville type relations for the long $h$-vectors}

We say, for brevity, that a face system
$\Phi\subset\mathbf{2}^{[m]}$ is a {\em DS-system} if the vector
$\pmb{h}(\Phi;m)$ satisfies the Dehn-Sommerville type relations
\begin{equation}
\label{eq:2} h_l(\Phi;m)=(-1)^{m-\size(\Phi)}h_{m-l}(\Phi;m)\ ,\ \
\ 0\leq l\leq m\ ;
\end{equation}
see also~\cite[\S{}7]{M1}.

If $\#\Phi>0$, then define the integer
\begin{equation*}
\eta(\Phi):=\begin{cases}|\bigcup_{F\in\Phi}F|, &\text{if
$|\bigcup_{F\in\Phi}F|\equiv \size(\Phi)\pmod{2}$,}\\
|\bigcup_{F\in\Phi}F|+1, &\text{if
$|\bigcup_{F\in\Phi}F|\not\equiv \size(\Phi)\pmod{2}$.}
\end{cases}
\end{equation*}
A face system $\Phi$ with $\#\Phi>0$ is a DS-system if and only if
for any $n\in\mathbb{P}$ such that
\begin{equation}
\label{eq:3}
\begin{split}
\eta(\Phi)&\leq n\ ,\\ n&\equiv\size(\Phi)\pmod{2}\ ,
\end{split}
\end{equation}
the vector $\pmb{h}(\Phi;n)$ is a left eigenvector (corresponding
to the eigenvalue $1$) of the $(n+1)\times(n+1)$ backward identity
matrix $\mathbf{U}(n)$:
\begin{equation*}
\pmb{h}(\Phi;n)=\pmb{h}(\Phi;n)\cdot\mathbf{U}(n)\ .
\end{equation*}

We denote the eigenspace of $\mathbf{U}(m)$ corresponding to the
eigenvalue~$1$ by $\boldsymbol{\mathcal{E}}^{\mathrm{h}}(m)$.

\subsection{The eigenvalues of a backward identity matrix}
\label{section:spectrum}

{\small Recall that the eigenvalues of the {\em permutation
matrix\/} $\mathbf{U}(m)$ are $-1$ and $1$; the matrix is {\em
unimodular}, see, e.g.,~\cite[Chapter~4]{Schrijver},
\cite[Chapter~8]{Sturmfels} on unimodular matrices. The algebraic
multiplicity of the eigenvalue $1$ is
$\left\lceil\tfrac{m+1}{2}\right\rceil$: the characteristic
polynomial $\wp\bigl(\mathbf{U}(m)\bigr)$ of the matrix, in the
variable $\lambda$, is
\begin{equation*}\wp\bigl(\mathbf{U}(m)\bigr)=
\begin{cases}(\lambda-1)^{\frac{m+2}{2}}
(\lambda+1)^{\frac{m}{2}},&\text{$m$ even},\\
(\lambda-1)^{\frac{m+1}{2}} (\lambda+1)^{\frac{m+1}{2}},&\text{$m$
odd}.
\end{cases}
\end{equation*}
In other words,
\begin{equation*}
\wp\bigl(\mathbf{U}(m)\bigr)=
\begin{cases}
\sum_{s=0}^{m+1} \mathbf{K}_s(m+1,\frac{m+2}{2})\cdot\lambda^s,
&\text{$m$ even},\\
\sum_{s=0}^{m+1}\mathbf{K}_s(m+1,\frac{m+1}{2})\cdot
\lambda^s,&\text{$m$ odd},
\end{cases}
\end{equation*}
where $\mathbf{K}_s(t,i)$ stands for the {\em Krawtchouk
polynomial\/} defined by $\sum_{s=0}^t\mathbf{K}_s(t,i)$
$\cdot\lambda^s:=(-\lambda+1)^i (\lambda+1)^{t-i}$; see,
e.g.,~\cite[\S{}1.2]{van Lint} on these polynomials.

The geometric multiplicity of $1$ equals its algebraic
multiplicity.
}

\subsection{The linear hulls of the long $h$-vectors of
DS-systems}

For a positive integer $k$, define the complex
$\overline{\mathbf{2}^{[k]}}:=\mathbf{2}^{[k]}-\{[k]\}
=[\hat{0},[k]]-\{[k]\}$.

\begin{proposition}
\begin{itemize}
\item[\rm(i)]
If $m$ is even, then
\begin{equation*}
\boldsymbol{\mathcal{E}}^{\mathrm{h}}(m)=
\lin\left(\pmb{h}(\overline{\mathbf{2}^{[1]}};m),\
\pmb{h}(\overline{\mathbf{2}^{[3]}};m),\ \ldots,\
\pmb{h}(\overline{\mathbf{2}^{[m-1]}};m),\
\boldsymbol{\iota}(m)\right)\ ;
\end{equation*}
if $m$ is odd, then
\begin{equation*}
\boldsymbol{\mathcal{E}}^{\mathrm{h}}(m)=
\lin\left(\pmb{h}(\overline{\mathbf{2}^{[2]}};m),\
\pmb{h}(\overline{\mathbf{2}^{[4]}};m),\ \ldots,\
\pmb{h}(\overline{\mathbf{2}^{[m-1]}};m),\
\boldsymbol{\iota}(m)\right)\ .
\end{equation*}

\item[\rm(ii)] If $\Phi\subset\mathbf{2}^{[m]}$ is a DS-system
such that $\Phi\not\ni[m]$ {\rm(}in particular, if $\Phi$ is a
complex whose $h$-vector satisfies {\rm(\ref{eq:12})\rm)} then
$\pmb{h}(\Phi;m)$ $\cdot\boldsymbol{\iota}(m)^{\top}=0$.
\end{itemize}
\end{proposition}

\begin{proof}
The proof of assertion~(i) is straightforward, with respect to the
argument from \S\ref{section:spectrum}. Assertion~(ii) is a
consequence of~\cite[Eq.~(3.8)]{M1}.
\end{proof}

Define a linear hyperplane $\boldsymbol{\mathcal{H}}(m)$ in the
space $\boldsymbol{\mathcal{E}}^{\mathrm h}(m)$, in the following
way: if $m$ is even, then
\begin{equation*}
\boldsymbol{\mathcal{H}}(m):=
\lin\left(\pmb{h}(\overline{\mathbf{2}^{[1]}};m),\
\pmb{h}(\overline{\mathbf{2}^{[3]}};m),\ \ldots,\
\pmb{h}(\overline{\mathbf{2}^{[m-1]}};m)\right)\ ,
\end{equation*}
with $\dim\boldsymbol{\mathcal{H}}(m)=\tfrac{m}{2}$; if $m$ is
odd, then
\begin{equation*}
\boldsymbol{\mathcal{H}}(m):=
\lin\left(\pmb{h}(\overline{\mathbf{2}^{[2]}};m),\
\pmb{h}(\overline{\mathbf{2}^{[4]}};m),\ \ldots,\
\pmb{h}(\overline{\mathbf{2}^{[m-1]}};m)\right)\ ,
\end{equation*}
$\dim\boldsymbol{\mathcal{H}}(m)=\tfrac{m-1}{2}$.

The one-dimensional subspace $\lin(\boldsymbol{\iota}(m))$ is the
orthogonal complement of the subspace
$\boldsymbol{\mathcal{H}}(m)$ of the space
$\boldsymbol{\mathcal{E}}^{\mathrm{h}}(m)$, with respect to the
standard scalar product.

\subsection{Some bases of $\mathbb{R}^{m+1}$}
Let $\{F_0,\ldots, F_m\}\subset\mathbf{2}^{[m]}$ be a face system
such that $|F_k|=k$, for $0\leq k\leq m$; here, $F_0:=\hat{0}$ and
$F_m:=[m]$. In~\cite[\S{}4]{M1} the following six bases of
$\mathbb{R}^{m+1}$ were considered:{\footnotesize
\begin{align*}
\bS_m=\bigl(\boldsymbol{\sigma}(0;m),\ldots,\boldsymbol{\sigma}(m;m)\bigr)
:&=\bigl(\ \pmb{f}(\{F_0\};m),\ldots,\pmb{f}(\{F_m\};m)\ \bigr)\
,\\ \bH_m^{\bullet}=\bigl(\boldsymbol{\vartheta}^{\bullet}(0;m),
\ldots,\boldsymbol{\vartheta}^{\bullet}(m;m)\bigr) :&=\bigl(\
\pmb{h}(\{F_0\};m),\ldots,\pmb{h}(\{F_m\};m)\ \bigr)\ ,\\
\bF_m^{\blacktriangle}=\bigl(\boldsymbol{\varphi}^{\blacktriangle}(0;m),\ldots,\boldsymbol{\varphi}^{\blacktriangle}(m;m)\bigr)
:&=\bigl(\ \pmb{f}([F_0,F_0];m),\ldots,\pmb{f}([F_0,F_m];m)\
\bigr)\ ,\\
\bH_m^{\blacktriangle}=\bigl(\boldsymbol{\vartheta}^{\blacktriangle}(0;m),\ldots,\boldsymbol{\vartheta}^{\blacktriangle}(m;m)\bigr)
:&=\bigl(\ \pmb{h}([F_0,F_0];m),\ldots,\pmb{h}([F_0,F_m];m)\
\bigr)\ ,\\
\bF_m^{\blacktriangledown}=\bigl(\boldsymbol{\varphi}^{\blacktriangledown}(0;m),\ldots,\boldsymbol{\varphi}^{\blacktriangledown}(m;m)\bigr)
:&=\bigl(\ \pmb{f}([F_m,F_m];m),\ldots,\pmb{f}([F_0,F_m];m)\
\bigr)\ ,\\
\bH_m^{\blacktriangledown}=\bigl(\boldsymbol{\vartheta}^{\blacktriangledown}(0;m),\ldots,\boldsymbol{\vartheta}^{\blacktriangledown}(m;m)\bigr)
:&=\bigl(\ \pmb{h}([F_m,F_m];m),\ldots,\pmb{h}([F_0,F_m];m)\
\bigr)\ ;
\end{align*}
} the basis $\bH_m^{\blacktriangledown}$ is up to rearrangement
the standard basis $\bS_m$.

Table~\ref{table:1} collects the representations of the vectors
$\pmb{h}(\overline{\mathbf{2}^{[k]}};m)$ with respect to the above
mentioned bases.

We define a matrix $\mathbf{S}(m)$ whose $(i,j)$th entry is
$(-1)^{j-i}\tbinom{m-i}{j-i}$ by {\footnotesize
\begin{equation*}
\mathbf{S}(m):=\begin{pmatrix}
\boldsymbol{\vartheta}^{\bullet}(0;m)\\ \vdots\\
\boldsymbol{\vartheta}^{\bullet}(m;m)
\end{pmatrix}\ ,
\end{equation*}
} and denote by $\mathrm{S}_m:\mathbb{R}^{m+1}\to\mathbb{R}^{m+1}$
the automorphism $\pmb{v}\mapsto\pmb{v}\cdot\mathbf{S}(m)$.

{\footnotesize
\begin{table}[ht]
\caption{Representations of
$\pmb{h}(\overline{\mathbf{2}^{[k]}};m)$ and
$\pmb{f}(\overline{\mathbf{2}^{[k]}};m)$, $1\leq k\leq m$, with
respect to various bases} \label{table:1}
\begin{center}
\begin{tabular}{|c|c|} \hline

\em $l$th component & \em Expression \\ \hline\hline

$h_l(\overline{\mathbf{2}^{[k]}};m)$ &
$\vartheta_l^{\blacktriangle}(k;m)-
\vartheta_l^{\bullet}(k;m)=(-1)^l\left(\binom{m-k}{l}
-(-1)^k\binom{m-k}{l-k}\right)$\\ \hline

$\kappa_l\bigl(\boldsymbol{\pmb{h}}(\overline{\mathbf{2}^{[k]}};m),
\bH^{\bullet}_m\bigr)$ & $\binom{k}{l}-\delta_{k,l}$\\ \hline

$\kappa_l\bigl(\boldsymbol{\pmb{h}}(\overline{\mathbf{2}^{[k]}};m),
\bF^{\blacktriangle}_m\bigr)$ & $(-1)^l\left(
2^{m-k-l}\binom{m-k}{l}-(-1)^k\sum_{s=l}^m\binom{m-k}{s-k}\binom{s}{l}\right)$\\
\hline

$\kappa_l\bigl(\boldsymbol{\pmb{h}}(\overline{\mathbf{2}^{[k]}};m),
\bH^{\blacktriangle}_m\bigr)$ &
$(-1)^{k-l}\left(\delta_{k,l}-\binom{k}{l}\right) $\\ \hline

$\kappa_l\bigl(\boldsymbol{\pmb{h}}(\overline{\mathbf{2}^{[k]}};m),
\bF^{\blacktriangledown}_m\bigr)$ & $(-1)^{m-l}\sum_{s=0}^m\left(
\binom{m-k}{s}-(-1)^k\binom{m-k}{s-k}\right)\binom{m-s}{l}$\\
\hline

$\kappa_l\bigl(\boldsymbol{\pmb{h}}(\overline{\mathbf{2}^{[k]}};m),
\bH^{\blacktriangledown}_m\bigr)$ &
$(-1)^{m-l}\left(\binom{m-k}{l-k}-(-1)^k\binom{m-k}{l}\right)$\\
\hline \hline


$f_l(\overline{\mathbf{2}^{[k]}};m)$ &
$\varphi_l^{\blacktriangle}(k;m)-
\sigma_l(k;m)=\varphi_l^{\blacktriangle}(k;m)-
\delta_{k,l}=\binom{k}{l} -\delta_{k,l}$\\ \hline

$\kappa_l\bigl(\boldsymbol{\pmb{f}}(\overline{\mathbf{2}^{[k]}};m),
\bH^{\bullet}_m\bigr)$ &
$\sum_{s=0}^{\min\{k,l\}}\binom{k}{s}\binom{m-s}{m-l}-\binom{m-k}{m-l}$\\
\hline

$\kappa_l\bigl(\boldsymbol{\pmb{f}}(\overline{\mathbf{2}^{[k]}};m),
\bF^{\blacktriangle}_m\bigr)$ & $(-1)^{k-l}\left(
\delta_{k,l}-\binom{k}{l}\right)$\\ \hline

$\kappa_l\bigl(\boldsymbol{\pmb{f}}(\overline{\mathbf{2}^{[k]}};m),
\bH^{\blacktriangle}_m\bigr)$ &
$(-1)^{m-l}\binom{k}{m-l}\left(2^{k+l-m}-1\right)$\\ \hline

$\kappa_l\bigl(\boldsymbol{\pmb{f}}(\overline{\mathbf{2}^{[k]}};m),
\bF^{\blacktriangledown}_m\bigr)$ &
$\binom{m-k}{m-l}-(-1)^{m-k-l}\binom{m-k}{l}$\\ \hline

$\kappa_l\bigl(\boldsymbol{\pmb{f}}(\overline{\mathbf{2}^{[k]}};m),
\bH^{\blacktriangledown}_m\bigr)$ &
$\binom{k}{m-l}-\delta_{k,m-l}$\\ \hline

\end{tabular}
\end{center}
\end{table}
}

\section{Dehn-Sommerville type relations for the long $f$-vectors}

Define a unimodular matrix $\mathbf{D}(m)$ by {\footnotesize
\begin{equation*}
\mathbf{D}(m):=\mathbf{S}(m)\cdot\mathbf{U}(m)\cdot\mathbf{S}(m)^{-1}=
\mathbf{S}(m)\cdot\begin{pmatrix}
\boldsymbol{\varphi}^{\blacktriangledown}(0;m)\\ \vdots \\
\boldsymbol{\varphi}^{\blacktriangledown}(m;m)\end{pmatrix}=
\begin{pmatrix}
\boldsymbol{\varphi}^{\blacktriangledown}(0;m)\\ \vdots \\
\boldsymbol{\varphi}^{\blacktriangledown}(m;m)
\end{pmatrix}^{-1}\cdot\mathbf{S}(m)^{-1}\ ;
\end{equation*}
} its $(i,j)$th entry is $(-1)^{m-i}\tbinom{i}{j}$.

Since the matrices $\mathbf{D}(m)$ and $\mathbf{U}(m)$ are
similar, the properties of $\mathbf{D}(m)$ coincide with those of
$\mathbf{U}(m)$ mentioned in \S\ref{section:spectrum}. Further,
for any DS-system $\Phi$ such that $\#\Phi>0$, and for any integer
$n$ satisfying~(\ref{eq:3}), it holds
\begin{equation}
\label{eq:4}\pmb{f}(\Phi;n)=\pmb{f}(\Phi;n)\cdot\mathbf{D}(n)\ ;
\end{equation}
in other words, we have
\begin{equation*}
f_l(\Phi;n)=(-1)^n\sum_{i=l}^{\size(\Phi)}(-1)^i\binom{i}{l}
f_i(\Phi;n)\ ,\ \ \ 0\leq l\leq\size(\Phi)\ ;
\end{equation*}
cf., e.g.,~\cite[p.~253]{Z}.

Let $\mathrm{S}_m^{-1}:\mathbb{R}^{m+1}\to\mathbb{R}^{m+1}$ be the
automorphism $\pmb{v}\mapsto\pmb{v}\cdot\mathbf{S}(m)^{-1}$. We
define subspaces
$\boldsymbol{\mathcal{F}}(m)\subset\boldsymbol{\mathcal{E}}^{\mathrm
f}(m)$ of $\mathbb{R}^{m+1}$ by $\boldsymbol{\mathcal{F}}(m):=
\mathrm{S}_m^{-1}\bigl(\boldsymbol{\mathcal{H}}(m)\bigr)$ and
$\boldsymbol{\mathcal{E}}^{\mathrm
f}(m):=\mathrm{S}_m^{-1}\bigl(\boldsymbol{\mathcal{E}}^{\mathrm
h}(m)\bigr)$. Thus, if $m$ is even, then
\begin{equation}
\label{eq:5}
\boldsymbol{\mathcal{F}}(m)=\lin\left(\pmb{f}(\overline{\mathbf{2}^{[1]}};m),\
\ \pmb{f}(\overline{\mathbf{2}^{[3]}};m),\ \ \ldots,\ \
\pmb{f}(\overline{\mathbf{2}^{[m-1]}};m)\right)\ ;
\end{equation}
if $m$ is odd, then
\begin{equation}
\label{eq:6}
\boldsymbol{\mathcal{F}}(m)=\lin\left(\pmb{f}(\overline{\mathbf{2}^{[2]}};m),\
\ \pmb{f}(\overline{\mathbf{2}^{[4]}};m),\ \ \ldots,\ \
\pmb{f}(\overline{\mathbf{2}^{[m-1]}};m)\right)\ .
\end{equation}

Define the row vector
\begin{equation*}
\boldsymbol{\pi}(m):=
(\tbinom{m+1}{0},\tbinom{m+1}{1},\ldots,\tbinom{m+1}{m}) =
\boldsymbol{\iota}(m)\cdot\mathbf{S}(m)^{-1}\in\mathbb{N}^{m+1}\ .
\end{equation*}

For any $m$, we have $\boldsymbol{\mathcal{E}}^{\mathrm
f}(m)=\boldsymbol{\mathcal{F}}(m)\oplus\lin(\boldsymbol{\pi}(m))$;
this is the eigenspace of the matrix $\mathbf{D}(m)$ corresponding
to its eigenvalue $1$.

All the vectors that appear in expressions~(\ref{eq:5})
and~(\ref{eq:6}), as well as the vector $\boldsymbol{\pi}(m)$, lie
in the affine hyperplane $\{\mathbf{z}\in\mathbb{R}^{m+1}:\
\mathrm{z}_0=1\}$. We also have
$(\boldsymbol{\pi}(m),0,0)=\pmb{f}(\overline{\mathbf{2}^{[m+1]}};m+2)$.
The representations of the vectors
$\pmb{f}(\overline{\mathbf{2}^{[k]}};m)$ with respect to various
bases are collected in Table~\ref{table:1}.

\section{The long $f$-vectors of DS-systems, and integer points in
rational polytopes}

Studying the long $f$-vectors of the DS-systems
$\Phi\subset\mathbf{2}^{[m]}$ such that $\#\Phi>0$ and
$m\equiv\size(\Phi)\pmod{2}$, and considering
relation~(\ref{eq:4}), we are interested in the solutions
$\mathbf{z}\in\mathbf{N}^{m+1}$, $\mathbf{z}\neq\pmb{0}$, to the
system
\begin{equation*}
\mathbf{z}\cdot\bigl(\mathbf{I}(m)-\mathbf{D}(m)\bigr)=\pmb{0}\ ,\
\ \ \pmb{0}\leq\mathbf{z}\leq(\tbinom{m}{0},
\tbinom{m}{1},\ldots,\tbinom{m}{m})\ .
\end{equation*}
The matrix $\mathbf{I}(m)-\mathbf{D}(m)$ is lower-triangular, of
rank $\left\lfloor\tfrac{m+1}{2}\right\rfloor$, with the
eigenvalues $0$ and $2$.

Define the polytope
\begin{multline*}
\boldsymbol{\mathcal{Q}}^{\mathrm{f}}(m):=\Bigl\{\mathbf{x}\in\mathbb{R}^{m+1}:\
\mathbf{x}\cdot\bigl(\mathbf{I}(m)-\mathbf{D}(m)\bigr)=\pmb{0}\
,\\ \pmb{0}\leq\mathbf{x}\leq(\tbinom{m}{0},\tbinom{m}{1},
\ldots,\tbinom{m}{m})\Bigr\}\ .
\end{multline*}

\begin{figure}
\begin{picture}(18,23)(-6,-17.5)

\put(5,-5){\makebox(0,0)[l]{\bf (a)}}

\put(0,0){\vector(0,1){5}} \put(0,0){\vector(3,-1){8}}
\put(0,0){\vector(-2,-3){3}}

\put(0,0){\line(3,2){3}}

\put(0,0){\line(5,2){7}} \put(7,2.8){\circle * {0.2}}

\put(0,0){\line(3,-4){3}} \put(3,-4){\circle * {0.2}}

\put(2.25,0.9){\circle * {0.1}}

\put(0,0){\circle * {0.2}} \put(3,2){\circle * {0.2}}
\put(1,-1){\circle * {0.2}}

\put(0,3){\circle * {0.15}} \put(6,-2){\circle * {0.1}}

\thicklines

\put(-2,-3){\line(0,1){3}} \put(-2,-3){\line(3,-1){6}}
\put(-2,0){\line(3,-1){6}} \put(4,-5){\line(0,1){3}}
\put(4,-5){\line(2,3){2}} \put(6,-2){\line(0,1){3}}
\put(0,3){\line(3,-1){6}} \put(-2,0){\line(2,3){2}}
\put(4,-2){\line(2,3){2}}

\put(-2,-3){\line(3,2){3}} \put(1,-1){\line(2,3){2}}

\put(-2,-3){\circle * {0.2}}

\put(-3,-4.75){\makebox(0,0)[r]{$\mathrm{x}_0$}}
\put(8.25,-2.75){\makebox(0,0)[l]{$\mathrm{x}_1$}}
\put(0,5.25){\makebox(0,0)[b]{$\mathrm{x}_2$}}

\put(7.25,2.75){\makebox(0,0)[l]{\scriptsize
$\boldsymbol{\pi}(2)=(1,3,3)$}}
\put(-0.25,0){\makebox(0,0)[r]{\scriptsize $(0,0,0)$}}
\put(-2.25,-3){\makebox(0,0)[r]{\scriptsize
$\pmb{f}(\overline{\mathbf{2}^{[1]}};2)=(1,0,0)$}}
\put(1.25,-0.9){\makebox(0,0)[l]{\scriptsize $(1,1,1)$}}
\put(3.25,2.25){\makebox(0,0)[l]{\scriptsize $(0,1,1)$}}

\put(0.25,3.25){\makebox(0,0)[l]{\scriptsize $(0,0,1)$}}
\put(6.25,-1.75){\makebox(0,0)[l]{\scriptsize $(0,2,0)$}}

\put(2.8,-3.9){\makebox(0,0)[r]{\scriptsize $(0,1,-1)$}}


\put(-5,-17){\makebox(0,0)[r]{\bf (b)}}

\thinlines

\put(0,-10){\vector(0,1){5}} \put(0,-10){\vector(3,-1){8}}
\put(0,-10){\vector(-2,-3){4}}

\put(0,-10){\line(3,-4){6}} \put(0,-10){\line(-4,1){8}}

\put(0,-10){\line(-1,-3){2}}

\put(0,-10){\line(1,-1){1}}

\thicklines

\put(6,-18){\line(0,1){3}}  \put(-2,-16){\line(0,1){3}}
\put(-8,-8){\line(0,1){3}}

\put(-8,-8){\line(3,-4){6}} \put(-2,-16){\line(4,-1){8}}
\put(-8,-5){\line(3,-4){6}} \put(-2,-13){\line(4,-1){8}}

\put(-8,-5){\line(4,-1){8}} \put(0,-7){\line(3,-4){6}}

\put(-4,-16.25){\makebox(0,0)[r]{$\mathrm{x}_0$}}
\put(8.25,-12.75){\makebox(0,0)[l]{$\mathrm{x}_1$}}
\put(0,-4.75){\makebox(0,0)[b]{$\mathrm{x}_2$}}

\put(3,-11){\line(-4,1){8}} \put(-8,-8){\line(3,-1){3}}

\put(0,-10){\circle * {0.2}} \put(-8,-8){\circle * {0.2}}
\put(-5,-9){\circle * {0.2}} \put(3,-11){\circle * {0.2}}

\put(-2,-16){\circle * {0.2}} \put(1,-11){\circle * {0.2}}

\put(0.35,-10.35){\circle * {0.1}}

\put(-7.75,-7.35){\makebox(0,0)[l]{\scriptsize
$\pmb{h}(\overline{\mathbf{2}^{[1]}};2)$}}

\put(-7.75,-7.75){\makebox(0,0)[l]{\scriptsize $=(1,-2,1)$}}

\put(-5.25,-9.25){\makebox(0,0)[r]{\scriptsize $(1,-1,1)$}}
\put(0.25,-9.75){\makebox(0,0)[l]{\scriptsize $(0,0,0)$}}
\put(3.25,-10.75){\makebox(0,0)[l]{\scriptsize $(0,1,0)$}}

\put(-2.25,-16.25){\makebox(0,0)[r]{\scriptsize $(1,0,-1)$}}

\put(1.25,-11){\makebox(0,0)[l]{\scriptsize
$\boldsymbol{\iota}(2)$}}
\put(1.25,-11.35){\makebox(0,0)[l]{\scriptsize $=(1,1,1)$}}

\end{picture}
\caption{$\quad$} \label{figure:1}
\end{figure}

If $m$ is even, then another description of
$\boldsymbol{\mathcal{Q}}^{\mathrm{f}}(m)$ is
\begin{multline*}
\boldsymbol{\mathcal{Q}}^{\mathrm{f}}(m)=
\lin\bigl(\pmb{f}(\overline{\mathbf{2}^{[1]}};m),\
\pmb{f}(\overline{\mathbf{2}^{[3]}};m),\ \ldots,\
\pmb{f}(\overline{\mathbf{2}^{[m-1]}};m),\ \boldsymbol{\pi}(m)
\bigr)
\\ \cap\ \ \left\{\mathbf{x}\in\mathbb{R}^{m+1}:\
\pmb{0}\leq\mathbf{x}\leq(\tbinom{m}{0},\tbinom{m}{1},
\ldots,\tbinom{m}{m})\right\}\ ;
\end{multline*}
if $m$ is odd, then we have
\begin{multline*}
\boldsymbol{\mathcal{Q}}^{\mathrm{f}}(m)=
\lin\bigl(\pmb{f}(\overline{\mathbf{2}^{[2]}};m),\
\pmb{f}(\overline{\mathbf{2}^{[4]}};m),\ \ldots,\
\pmb{f}(\overline{\mathbf{2}^{[m-1]}};m),\ \boldsymbol{\pi}(m)
\bigr)
\\ \cap\ \ \left\{\mathbf{x}\in\mathbb{R}^{m+1}:\
\pmb{0}\leq\mathbf{x}\leq(\tbinom{m}{0},\tbinom{m}{1},
\ldots,\tbinom{m}{m})\right\}\ ;
\end{multline*}
In other words, for any $m$, we have
\begin{equation*}
\boldsymbol{\mathcal{Q}}^{\mathrm{f}}(m) =
\boldsymbol{\mathcal{E}}^{\mathrm{f}}(m)\ \cap\
\boldsymbol{\varPi}(m)\ ,
\end{equation*}
where {\small
\begin{equation}
\label{eq:15}
\boldsymbol{\varPi}(m):=\conv\left(\sum_{k=0}^ma_k\binom{m}{k}
\boldsymbol{\sigma}(k;m):\
(a_0,\ldots,a_m)\in\{0,1\}^{m+1}\right)\ .
\end{equation}
}

On the one hand, the points $\{\mathbf{z}\in
\boldsymbol{\mathcal{Q}}^{\mathrm{f}}(m)\cap\mathbb{N}^{m+1}:\
\mathbf{z}\neq\pmb{0}\}$ are exactly the vectors $\pmb{f}(\Phi;m)$
of the DS-systems $\Phi\subset\mathbf{2}^{[m]}$ such that
$\#\Phi>0$ and $m\equiv\size(\Phi)\pmod{2}$. On the other hand, if
$\mathbf{z}'\in\boldsymbol{\mathcal{Q}}^{\mathrm{f}}(m)$ then
there are $\prod_{k=0}^m\binom{\tbinom{m}{k}}{\mathrm{z}'_k}$
DS-systems $\Phi\subset\mathbf{2}^{[m]}$ corresponding to the
point $\mathbf{z}'$.

\begin{example}{\small
We have
$\boldsymbol{\mathcal{Q}}^{\mathrm{f}}(2)=\conv\bigl(\pmb{0},(1,0,0),
(0,1,1),(1,1,1)\bigr)$, see Figure~{\rm\ref{figure:1}(a)}.
}
\end{example}

One more description of $\boldsymbol{\mathcal{Q}}^{\mathrm{f}}(m)$
is
\begin{equation*}
\boldsymbol{\mathcal{Q}}^{\mathrm{f}}(m)=\left\{\mathbf{x}\in
\boldsymbol{\mathcal{C}}^{\mathrm{f}}(m):\
\mathbf{x}\leq(\tbinom{m}{0},\tbinom{m}{1},
\ldots,\tbinom{m}{m})\right\}\ ,
\end{equation*}
where a convex pointed polyhedral cone
$\boldsymbol{\mathcal{C}}^{\mathrm{f}}(m)$ is defined by
\begin{equation*}
\boldsymbol{\mathcal{C}}^{\mathrm{f}}(m):=
\boldsymbol{\mathcal{E}}^{\mathrm f}(m)\ \cap\
\cone\bigl(\boldsymbol{\sigma}(0;m),\ \boldsymbol{\sigma}(1;m),\
\ldots,\ \boldsymbol{\sigma}(m;m)\bigr) \ .
\end{equation*}

\begin{proposition}
\label{prop:1}
\begin{itemize}
\item[\rm(i)]
If $m$ is even, then
\begin{equation*}
\boldsymbol{\mathcal{C}}^{\mathrm{f}}(m)\supseteq\cone
\left(\boldsymbol{\varphi}^{\blacktriangle}(i;m)\cdot\mathbf{T}(m)^i:\
0\leq i\leq\tfrac{m}{2}\right)\ .
\end{equation*}
For any $i\in\mathbb{N}$, $i\leq\tfrac{m}{2}$, the ray
$\cone\bigl(\boldsymbol{\varphi}^{\blacktriangle}(i;m)\cdot\mathbf{T}(m)^i\bigr)$
is an extreme ray of the $\tfrac{m+2}{2}$-dimensional unimodular
cone
\begin{equation}
\label{eq:14} \cone
\bigl(\boldsymbol{\varphi}^{\blacktriangle}(i;m)
\cdot\mathbf{T}(m)^i:\ 0\leq i\leq\tfrac{m}{2}\bigr)\ .
\end{equation}

\item[\rm(ii)] If $m$ is odd, then
\begin{multline*}
\boldsymbol{\mathcal{C}}^{\mathrm{f}}(m)\supseteq \cone
\Bigl(\boldsymbol{\varphi}^{\blacktriangle}(i+1;m)\cdot\mathbf{T}(m)^i\\+
\boldsymbol{\varphi}^{\blacktriangle}(i;m)\cdot\mathbf{T}(m)^{i+1}:\
0\leq i\leq\tfrac{m-1}{2}\Bigr)\ .
\end{multline*}
For any $i\in\mathbb{N}$, $i\leq\tfrac{m-1}{2}$, the ray $\cone
\bigl(\boldsymbol{\varphi}^{\blacktriangle}(i+1;m)\cdot\mathbf{T}(m)^i+
\boldsymbol{\varphi}^{\blacktriangle}(i;m)\cdot\mathbf{T}(m)^{i+1}\bigr)$
is an extreme ray of the $\tfrac{m+1}{2}$-dimensional unimodular
cone $\cone
\bigl(\boldsymbol{\varphi}^{\blacktriangle}(i+1;m)\cdot\mathbf{T}(m)^i+
\boldsymbol{\varphi}^{\blacktriangle}(i;m)\cdot\mathbf{T}(m)^{i+1}:\
0\leq i\leq\tfrac{m-1}{2}\bigr)$.
\end{itemize}
\end{proposition}

\begin{proof}
(i) The vectors from the sequence
\begin{equation}
\label{eq:13}
\bigl(\boldsymbol{\varphi}^{\blacktriangle}(i;m)\cdot\mathbf{T}(m)^i:\
0\leq i\leq\tfrac{m}{2}\bigr)
\end{equation}
are linearly independent, and for $i,j\in\mathbb{N}$ such that
$i\leq\tfrac{m}{2}$ and $j\leq m$, we have
\begin{equation}
\label{eq:7} \underbrace{\kappa_j}_{\text{or}\
\kappa_{m-j}}\left(\boldsymbol{\varphi}^{\blacktriangle}(i;m)\cdot\mathbf{T}(m)^i
\cdot\mathbf{S}(m)\right)=(-1)^{j-i}\tbinom{m-2i}{j-i}\ ,
\end{equation}
hence sequence~(\ref{eq:13}) is a basis of the space
$\boldsymbol{\mathcal{E}}^{\mathrm f}(m)$, and the cone generated
by it is simplicial. The matrix
\begin{equation*}
\begin{pmatrix}
\kappa_0\bigl(\boldsymbol{\varphi}^{\blacktriangle}(0;m)\cdot\mathbf{T}(m)^0\bigr)&
\cdots&\kappa_{\frac{m}{2}}\bigl(\boldsymbol{\varphi}^{\blacktriangle}(0;m)\cdot
\mathbf{T}(m)^0\bigr)\\ \hdotsfor{3}\\
\kappa_0\bigl(\boldsymbol{\varphi}^{\blacktriangle}(\tfrac{m}{2};m)\cdot
\mathbf{T}(m)^{\frac{m}{2}}\bigr)&
\cdots&\kappa_{\frac{m}{2}}\bigl(\boldsymbol{\varphi}^{\blacktriangle}(\tfrac{m}{2};m)
\cdot\mathbf{T}(m)^{\frac{m}{2}}\bigr)
\end{pmatrix}
\end{equation*}
is upper-triangular whose diagonal entries are $1$. This implies
that~(\ref{eq:13}) is an integral basis of the intersection of the
linear hull of~(\ref{eq:13}) with $\mathbb{Z}^{m+1}$. In other
words, cone~(\ref{eq:14}) is unimodular.

Pick a vector $\pmb{v}\in\boldsymbol{\mathcal{E}}^{\mathrm f}(m)$
such that $\pmb{v}\geq\pmb{0}$ and $\pmb{v}\neq\pmb{0}$. If
$\pmb{v}=\sum_{i=0}^{m/2}a_i$
$\cdot\boldsymbol{\varphi}^{\blacktriangle}(i;m)
\cdot\mathbf{T}(m)^i$, for some
$a_0,\ldots,a_{\frac{m}{2}}\in\mathbb{R}$, then the equalities
$\kappa_0(\pmb{v})=a_0$, $\kappa_1(\pmb{v})=a_1$ and
$\kappa_{\frac{m}{2}}(\pmb{v})=a_{\frac{m}{2}}$ imply that
$a_0,a_1,a_{\frac{m}{2}}\geq 0$. As a consequence, if
$m\in\{2,4\}$ then
\begin{equation*}
\boldsymbol{\mathcal{C}}^{\mathrm{f}}(m)=\cone
\left(\boldsymbol{\varphi}^{\blacktriangle}(i;m)\cdot\mathbf{T}(m)^i:\
0\leq i\leq\tfrac{m}{2}\right)\ .
\end{equation*}

The proof of assertion~(ii) is analogous to that of~(i). In
particular, for $i,j\in\mathbb{N}$ such that $i\leq\tfrac{m-1}{2}$
and $j\leq m$, we have
\begin{multline}
\label{eq:8} \underbrace{\kappa_j}_{\text{or}\
\kappa_{m-j}}\Bigl(\
\bigl(\boldsymbol{\varphi}^{\blacktriangle}(i+1;m)\cdot\mathbf{T}(m)^i+
\boldsymbol{\varphi}^{\blacktriangle}(i;m)\cdot\mathbf{T}(m)^{i+1}\bigr)
\cdot\mathbf{S}(m)\ \Bigr)\\=
(-1)^{j-i}(\tbinom{m-2i-1}{j-i}-\tbinom{m-2i-1}{j-i-1})\ ;
\end{multline}
the sequence
\begin{equation}
\label{eq:16}
\Bigl(\boldsymbol{\varphi}^{\blacktriangle}(i+1;m)\cdot\mathbf{T}(m)^i\\+
\boldsymbol{\varphi}^{\blacktriangle}(i;m)\cdot\mathbf{T}(m)^{i+1}:\
0\leq i\leq\tfrac{m-1}{2}\Bigr)
\end{equation}
is a basis of $\boldsymbol{\mathcal{E}}^{\mathrm f}(m)$.

If
$\pmb{v}:=\sum_{i=0}^{(m-1)/2}a_i\cdot\bigl(\boldsymbol{\varphi}^{\blacktriangle}(i+1;m)\cdot\mathbf{T}(m)^i+
\boldsymbol{\varphi}^{\blacktriangle}(i;m)\cdot\mathbf{T}(m)^{i+1}\bigr)\in
\boldsymbol{\mathcal{E}}^{\mathrm f}(m)$ and $\pmb{v}\geq\pmb{0}$
and $\pmb{v}\neq\pmb{0}$, for some
$a_0,\ldots,a_{\frac{m-1}{2}}\in\mathbb{R}$, then the equalities
$\kappa_0(\pmb{v})=a_0$ and
$\kappa_{\frac{m-1}{2}}(\pmb{v})=a_{\frac{m-1}{2}}$ imply that
$a_0,a_{\frac{m-1}{2}}\geq 0$ and, as a consequence, it holds
\begin{equation*}
\boldsymbol{\mathcal{C}}^{\mathrm{f}}(3)=\cone
\Bigl(\boldsymbol{\varphi}^{\blacktriangle}(i+1;3)\cdot\mathbf{T}(3)^i+
\boldsymbol{\varphi}^{\blacktriangle}(i;3)\cdot\mathbf{T}(3)^{i+1}:\
0\leq i\leq 1\Bigr)\ .
\end{equation*}
\end{proof}

The structure of bases~(\ref{eq:13}) and~(\ref{eq:16}) of
$\boldsymbol{\mathcal{E}}^{\mathrm f}(m)$ allows us to come to the
following conclusion:

\begin{corollary} \label{cor:1} $\quad$
\begin{itemize}
\item[\rm(i)]
Let $\Phi\subset\mathbf{2}^{[m]}$ be a DS-system.

If $m$ is even, then either $f_{m-1}(\Phi;m)=\tfrac{m}{2}$ and
$f_m(\Phi;m)=1$, or $f_{m-1}(\Phi;m)=f_m(\Phi;m)=0$; if $m$ is
odd, then it holds $f_{m-1}(\Phi;m)=f_m(\Phi;m)=0$.

\item[\rm(ii)]
If $m$ is even, then for any $t\in\mathbb{N}$,
$t\leq\tfrac{m-2}{2}$, it holds
\begin{equation*}
\lin\left(\pmb{f}(\overline{\mathbf{2}^{[2i+1]}};m):\ 0\leq i\leq
t\right)=\lin\left(\boldsymbol{\varphi}^{\blacktriangle}(i;m)
\cdot\mathbf{T}(m)^i:\ 0\leq i\leq t\right)\ ;
\end{equation*}
if $m$ is odd, then for any $t\in\mathbb{N}$,
$t\leq\tfrac{m-3}{2}$, it holds
\begin{multline*}
\lin\left(\pmb{f}(\overline{\mathbf{2}^{[2(i+1)]}};m):\ 0\leq
i\leq t\right)=
\lin\Bigl(\boldsymbol{\varphi}^{\blacktriangle}(i+1;m)\cdot\mathbf{T}(m)^i\\+
\boldsymbol{\varphi}^{\blacktriangle}(i;m)\cdot\mathbf{T}(m)^{i+1}:\
0\leq i\leq t\Bigr)\ .
\end{multline*}

\item[\rm(iii)]
{\small
\begin{equation*}
\lin\bigl(\pmb{f}(\Phi;m):\ \Phi\subset\mathbf{2}^{[m]},\
m\equiv\size(\Phi)\pmod{2},\ \Phi\text{\ is DS}\bigr)
=\begin{cases} \boldsymbol{\mathcal{E}}^{\mathrm{f}}(m),
&\text{$m$ even},\\ \boldsymbol{\mathcal{F}}(m), &\text{$m$ odd}.
\end{cases}
\end{equation*}
}
\end{itemize}
\end{corollary}

{\footnotesize
\begin{table}[ht]
\caption{Representations of
$\pmb{w}:=\boldsymbol{\varphi}^{\blacktriangle}(i;m)\cdot\mathbf{T}(m)^i$,
where $m$ is even and $0\leq i\leq \tfrac{m}{2}$, with respect to
various bases } \label{table:2}
\begin{center}
\begin{tabular}{|c|c|} \hline

\em $l$th component & \em Expression \\ \hline\hline

$\kappa_l\bigl(\pmb{w}, \bS_m\bigr)$ & $\binom{i}{l-i}$\\ \hline

$\kappa_l\bigl(\pmb{w}, \bH^{\bullet}_m\bigr)$ &
$\sum_{s=i}^{\min\{2i,l\}}\binom{i}{s-i}\binom{m-s}{m-l}$\\ \hline

$\kappa_l\bigl(\pmb{w}, \bF^{\blacktriangle}_m\bigr)$ &
$(-1)^l\binom{i}{l-i}$\\ \hline

$\kappa_l\bigl(\pmb{w}, \bH^{\blacktriangle}_m\bigr)$ &
$(-1)^{l}\sum_{s=\max\{i,m-l\}}^{2i}\binom{i}{s-i}
\binom{s}{m-l}$\\ \hline

$\kappa_l\bigl(\pmb{w}, \bF^{\blacktriangledown}_m\bigr)$ &
$(-1)^{l-i}\binom{m-2i}{l-i}$\\ \hline

$\kappa_l\bigl(\pmb{w}, \bH^{\blacktriangledown}_m\bigr)$ &
$\binom{i}{m-l-i}$\\ \hline \hline


$\kappa_l\bigl(\pmb{w}\cdot\mathbf{S}(m), \bS_m\bigr)$ &
$(-1)^{l-i}\binom{m-2i}{l-i}$\\ \hline

$\kappa_l\bigl(\pmb{w} \cdot\mathbf{S}(m), \bH^{\bullet}_m\bigr)$
& $\binom{i}{l-i}$\\ \hline

$\kappa_l\bigl(\pmb{w} \cdot\mathbf{S}(m),
\bF^{\blacktriangle}_m\bigr)$ &
$(-1)^{l-i}\sum_{s=\max\{i,l\}}^{m-i}\binom{s}{l}
\binom{m-2i}{s-i}$\\ \hline

$\kappa_l\bigl(\pmb{w} \cdot\mathbf{S}(m),
\bH^{\blacktriangle}_m\bigr)$ & $(-1)^l \binom{i}{l-i}$\\ \hline

$\kappa_l\bigl(\pmb{w} \cdot\mathbf{S}(m),
\bF^{\blacktriangledown}_m\bigr)$ &
$(-1)^{l-i}\sum_{s=\max\{i,l\}}^{m-i}\binom{s}{l}
\binom{m-2i}{s-i}$\\ \hline

$\kappa_l\bigl(\pmb{w} \cdot\mathbf{S}(m),
\bH^{\blacktriangledown}_m\bigr)$ &
$(-1)^{l-i}\binom{m-2i}{l-i}$\\ \hline

\end{tabular}
\end{center}
\end{table}
}

{\footnotesize
\begin{table}[ht]
\caption{Representations of
$\pmb{w}:=\boldsymbol{\varphi}^{\blacktriangle}(i+1;m)\cdot\mathbf{T}(m)^i+
\boldsymbol{\varphi}^{\blacktriangle}(i;m)\cdot\mathbf{T}(m)^{i+1}$,
where $m$ is odd and $0\leq i\leq\tfrac{m-1}{2}$, with respect to
various bases} \label{table:3}
\begin{center}
\begin{tabular}{|c|c|} \hline

\em $l$th component & \em Expression \\ \hline\hline

$\kappa_l\bigl(\pmb{w}, \bS_m\bigr)$ &
$\binom{i+1}{l-i}+\binom{i}{l-i-1}$\\ \hline

$\kappa_l\bigl(\pmb{w}, \bH^{\bullet}_m\bigr)$ & $\binom{m-i}{m-l}
+\sum_{s=i+1}^{\min\{2i+1,l\}}(\binom{i+1}{s-i}+\binom{i}{s-i-1})
\binom{m-s}{m-l}$\\ \hline

$\kappa_l\bigl(\pmb{w}, \bF^{\blacktriangle}_m\bigr)$ &
$(-1)^{l-1}(\binom{i}{l-i-1}+\binom{i+1}{l-i})$
\\ \hline

$\kappa_l\bigl(\pmb{w}, \bH^{\blacktriangle}_m\bigr)$ &
$(-1)^{l-1} \sum_{s=\max\{i,m-l\}}^{2i+1}(\binom{i+1}{s-i}+
\binom{i}{s-i-1})\tbinom{s}{m-l}$\\ \hline

$\kappa_l\bigl(\pmb{w}, \bF^{\blacktriangledown}_m\bigr)$ &
$(-1)^{l-i-1}(\tbinom{m-2i-1}{m-l-i}-\tbinom{m-2i-1}{m-l-i-1})$\\
\hline

$\kappa_l\bigl(\pmb{w}, \bH^{\blacktriangledown}_m\bigr)$ &
$\binom{i+1}{m-l-i}+\binom{i}{m-l-i-1}$\\ \hline \hline


$\kappa_l\bigl(\pmb{w}\cdot\mathbf{S}(m), \bS_m\bigr)$ &
$(-1)^{l-i-1}(\tbinom{m-2i-1}{m-l-i}-\tbinom{m-2i-1}{m-l-i-1})$\\
\hline

$\kappa_l\bigl(\pmb{w}\cdot\mathbf{S}(m), \bH^{\bullet}_m\bigr)$ &
$\binom{i}{l-i-1}+\binom{i+1}{l-i}$\\ \hline

$\kappa_l\bigl(\pmb{w}\cdot\mathbf{S}(m),
\bF^{\blacktriangle}_m\bigr)$ & $(-1)^{l-i}
\sum_{s=\max\{i,l\}}^{m-i}(\binom{m-2i-1}{s-i}-\binom{m-2i-1}{s-i-1})
\binom{s}{l} $\\ \hline

$\kappa_l\bigl(\pmb{w}\cdot\mathbf{S}(m),
\bH^{\blacktriangle}_m\bigr)$ &
$(-1)^{l-1}(\binom{i}{l-i-1}+\binom{i+1}{l-i})$\\ \hline

$\kappa_l\bigl(\pmb{w}\cdot\mathbf{S}(m),
\bF^{\blacktriangledown}_m\bigr)$ & $(-1)^{l-i}
\sum_{s=\max\{i,l\}}^{m-i}(\binom{m-2i-1}{s-i}-\binom{m-2i-1}{s-i-1})
\binom{s}{l} $\\ \hline

$\kappa_l\bigl(\pmb{w}\cdot\mathbf{S}(m),
\bH^{\blacktriangledown}_m\bigr)$ &
$(-1)^{l-i-1}(\tbinom{m-2i-1}{m-l-i}-\tbinom{m-2i-1}{m-l-i-1})$\\
\hline
\end{tabular}
\end{center}
\end{table}
}

If $m$ is even, define the $\tfrac{m}{2}$-dimensional polytope
\begin{multline*}
\boldsymbol{\mathcal{P}}^{\mathrm{f}}(m):=\boldsymbol{\mathcal{E}}^{\mathrm{f}}(m)\cap
\bigl\{\mathbf{z}\in\mathbb{R}^{m+1}:\ z_0=0,\\
\pmb{0}\leq(z_1,\ldots,z_m)\leq(\tbinom{m}{1},\tbinom{m}{2},\ldots,\tbinom{m}{m})
\bigr\}\ ,
\end{multline*}
and note that the polytope
$\boldsymbol{\mathcal{Q}}^{\mathrm{f}}(m)$ is a prism whose basis
is $\boldsymbol{\mathcal{P}}^{\mathrm{f}}(m)$:
\begin{equation*}
\boldsymbol{\mathcal{Q}}^{\mathrm{f}}(m)=
\boldsymbol{\mathcal{P}}^{\mathrm{f}}(m)\boxplus
\boldsymbol{\varphi}^{\blacktriangle}(0;m)\cdot\mathbf{T}(m)^0 =
\boldsymbol{\mathcal{P}}^{\mathrm{f}}(m)\boxplus
\boldsymbol{\sigma}(0;m)\ ;
\end{equation*}
as a consequence, we have
\begin{equation*}
\sum_{\boldsymbol{\alpha}\in\boldsymbol{\mathcal{Q}}^{\mathrm{f}}(m)\cap\mathbb{N}^{m+1}}
\mathbf{x}^{\boldsymbol{\alpha}}=(1+\mathrm{x}_0)\sum_{\boldsymbol{\alpha}\in\boldsymbol{\mathcal{P}}^{\mathrm{f}}(m)\cap\mathbb{N}^{m+1}}
\mathbf{x}^{\boldsymbol{\alpha}}\ ,
\end{equation*}
where
$\mathbf{x}^{\boldsymbol{\alpha}}:=\mathrm{x}_0{}^{\alpha_0}\cdots
\mathrm{x}_m{}^{\alpha_m}$.

Since for any $m$, the generating function for the long
$f$-vectors of DS-systems contained in $\mathbf{2}^{[m]}$ is
\begin{equation*}
\sum_{\substack{\boldsymbol{\mathfrak{f}}\in\mathbb{N}^{m+1}:\\
\exists\Phi\subset\mathbf{2}^{[m]},\\ \Phi \text{\ is DS,\ }
\pmb{f}(\Phi;m)=\boldsymbol{\mathfrak{f}}}}
\mathbf{x}^{\boldsymbol{\mathfrak{f}}}=-1+
\sum_{\boldsymbol{\alpha}\in
\boldsymbol{\mathcal{Q}}^{\mathrm{f}}(m)\cap\mathbb{N}^{m+1}}\mathbf{x}^{\boldsymbol{\alpha}}
+ \sum_{\substack{(\boldsymbol{\alpha},0)\in
\boldsymbol{\mathcal{Q}}^{\mathrm{f}}(m+1)\cap\mathbb{N}^{m+2}:\\
\boldsymbol{\alpha}\leq(\binom{m}{0},\binom{m}{1},\ldots,\binom{m}{m})
}}\mathbf{x}^{\boldsymbol{\alpha}}\ ,
\end{equation*}
we come to the conclusion:
\begin{proposition}
\begin{itemize}
\item[\rm(i)] If $m$ is even, then
\begin{equation*}
\begin{split}
\sum_{\substack{\boldsymbol{\mathfrak{f}}\in\mathbb{N}^{m+1}:\\
\exists\Phi\subset\mathbf{2}^{[m]},\\ \Phi \text{\rm\ is DS,\ }
\pmb{f}(\Phi;m)=\boldsymbol{\mathfrak{f}}}}
\mathbf{x}^{\boldsymbol{\mathfrak{f}}}=
-1&+(1+\mathrm{x}_0)\sum_{(0,\boldsymbol{\beta})\in
\boldsymbol{\mathcal{P}}^{\mathrm{f}}(m)\cap\mathbb{N}^{m+1}}
\mathbf{x}^{\boldsymbol{\beta}}\\&+
\sum_{\substack{(\boldsymbol{\gamma},0,0)\in
\boldsymbol{\mathcal{Q}}^{\mathrm{f}}(m+1)\cap\mathbb{N}^{m+2}:\\
\boldsymbol{\gamma}\leq(\binom{m}{0},\binom{m}{1},\ldots,\binom{m}{m-1})
}} \mathbf{x}^{\boldsymbol{\gamma}}\ .
\end{split}
\end{equation*}

\item[\rm(ii)] If $m$ is odd, then
\begin{equation*}
\begin{split}
\sum_{\substack{\boldsymbol{\mathfrak{f}}\in\mathbb{N}^{m+1}:\\
\exists\Phi\subset\mathbf{2}^{[m]},\\ \Phi \text{\rm\ is DS,\ }
\pmb{f}(\Phi;m)=\boldsymbol{\mathfrak{f}}}}
\mathbf{x}^{\boldsymbol{\mathfrak{f}}}
=-1&+\sum_{(\boldsymbol{\beta},0,0)\in
\boldsymbol{\mathcal{Q}}^{\mathrm{f}}(m)\cap\mathbb{N}^{m+1}}
\mathbf{x}^{\boldsymbol{\beta}}\\&+
(1+\mathrm{x}_0)\sum_{\substack{(0,\boldsymbol{\gamma},0,0)\in
\boldsymbol{\mathcal{P}}^{\mathrm{f}}(m+1)\cap\mathbb{N}^{m+2}:\\
\boldsymbol{\gamma}\leq(\binom{m}{1},\binom{m}{2},\ldots,\binom{m}{m-1})
}} \mathbf{x}^{\boldsymbol{\gamma}}\ .
\end{split}
\end{equation*}
\end{itemize}
\end{proposition}

{\footnotesize
\begin{table}[ht]
\caption{The number of the long $f$-vectors of the DS-systems
contained in $\mathbf{2}^{[m]}$, $2\leq m\leq 10$} \label{table:4}
\begin{center}
\begin{tabular}{|c|c|c|c|} \hline

$\quad$ & $|\{\ \boldsymbol{\mathfrak{f}}\in\mathbb{N}^{m+1}:$ &
$|\{\ \boldsymbol{\mathfrak{f}}\in\mathbb{N}^{m+1}:$ & $|\{\
\boldsymbol{\mathfrak{f}}\in\mathbb{N}^{m+1}:$ \\ $\quad$ &
$\exists\Phi\subset\mathbf{2}^{[m]}$, $\#\Phi>0$, &
$\exists\Phi\subset\mathbf{2}^{[m]}$, $\#\Phi>0$, &
$\exists\Phi\subset\mathbf{2}^{[m]}$, \\ $m$ &
$m\equiv\size(\Phi)\pmod{2}$, $\Phi$ is DS, &
$m\equiv\size(\Phi)+1\pmod{2}$, $\Phi$ is DS, & $\Phi$ is DS,
\\ $\quad$ & $\pmb{f}(\Phi;m)=\boldsymbol{\mathfrak{f}} \ \} |$ &
$\pmb{f}(\Phi;m)=\boldsymbol{\mathfrak{f}} \ \} |$ &
$\pmb{f}(\Phi;m)=\boldsymbol{\mathfrak{f}} \ \} |$
\\ \hline\hline
$2$ & $3$ & $1$ & $5$\\ $3$ & $1$ & $7$ & $9$ \\ $4$ & $19$ & $5$
& $25$\\ $5$ & $7$ & $71$ & $79$ \\ $6$ & $291$ & $41$ & $333$ \\
$7$ & $103$ & $2223$ & $2327$\\ $8$ & $17465$ & $1107$ & $18573$\\
$9$ & $4905$ & $271619$ & $276525$\\ $10$ & $3959091$ & $103057$ &
$4062149$\\ \hline
\end{tabular}
\end{center}
\end{table}
}

See, e.g., \cite[Chapter~VIII]{BarvinokConv}, \cite{BaPo},
\cite{BW}, \cite{BR}, \cite{BZ}, \cite{Brion}, \cite{LZ},
\cite[Chapter~12]{MS} on lattice-point counting in polytopes.
Effective tools for manipulating rational generating functions are
presented in~\cite{BW}.

Table~\ref{table:4} collects some illustrative information
obtained with the help of the software {\tt
LattE}~\cite{ShortLattE}, \cite{LattE}, \cite{EffectiveLattE}.


\section{Appendices}

{\small

\subsection{The long $h$-vectors of DS-systems,
and integer points in rational polytopes}

We now consider the polytope $\boldsymbol{\mathcal{Q}}^{\mathrm
h}(m):=\mathrm{S}_m\bigl(\boldsymbol{\mathcal{Q}}^{\mathrm
f}(m)\bigr):=\{\mathbf{z}\cdot\mathbf{S}(m):\
\mathbf{z}\in\boldsymbol{\mathcal{Q}}^{\mathrm f}(m)\}$.

The relation $\pmb{h}(\Phi;m)=\pmb{h}(\Phi;m)\cdot\mathbf{U}(m)$,
$\pmb{h}(\Phi;m)\neq\pmb{0}$, describes the long $h$-vectors of
the DS-systems $\Phi\subset\mathbf{2}^{[m]}$ such that $\#\Phi>0$
and $m\equiv\size(\Phi)\pmod{2}$, and  we are interested in the
solutions $\mathbf{z}\in\mathbb{Z}^{m+1}$,
$\mathbf{z}\neq\pmb{0}$, to the system
\begin{equation}
\label{eq:9}
\begin{split}
\mathbf{z}\cdot\bigl(\mathbf{I}(m)-\mathbf{U}(m)\bigr)=\pmb{0}\ ,
\\ \mathbf{z}\in\mathrm{S}_m\bigl(\boldsymbol{\varPi}(m)\bigr)=
\conv\left(\sum_{k=0}^ma_k\binom{m}{k}
\boldsymbol{\vartheta}^{\bullet}(k;m):\
(a_0,\ldots,a_m)\in\{0,1\}^{m+1}\right)\ .
\end{split}
\end{equation}

The matrix $\mathbf{I}(m)-\mathbf{U}(m)$ of rank
$\left\lfloor\tfrac{m+1}{2}\right\rfloor$ is {\em totally
unimodular} (see, e.g.,~\cite[\S{}19]{Schrijver} on totally
unimodular matrices,) that is, any its minor is either $-1$ or $0$
or $1$; this matrix can be substituted in system~(\ref{eq:9}) by
the submatrix composed of the first
$\left\lfloor\tfrac{m+1}{2}\right\rfloor$ columns.

System~(\ref{eq:9}) describes the intersection
\begin{equation}
\label{eq:10}
\boldsymbol{\mathcal{Q}}^{\mathrm{h}}(m):=\bigcap_{1\leq
k\leq\lfloor(m+1)/2\rfloor} \{\mathbf{x}\in\mathbb{R}^{m+1}:\
\mathrm{x}_{k-1}=\mathrm{x}_{m-k+1}\}\ \cap\
\mathrm{S}_m\bigl(\boldsymbol{\varPi}(m)\bigr)
\end{equation}
of the $\left\lceil\tfrac{m+1}{2}\right\rceil$-dimensional {\em
center\/} of a {\em graphical hyperplane arrangement} (see,
e.g.,~\cite[\S{}2.4]{OT}, \cite[Lecture~2]{St3} on graphical
arrangements) with an $(m+1)$-dimensional parallelepiped. The {\em
intersection poset\/} of the arrangement
$\bigl\{\{\mathbf{x}\in\mathbb{R}^{m+1}:\
\mathrm{x}_{k-1}=\mathrm{x}_{m-k+1}\}:\ 1\leq k\leq
\left\lfloor\tfrac{m+1}{2}\right\rfloor\bigr\}$ is a Boolean
lattice of rank $\left\lfloor\tfrac{m+1}{2}\right\rfloor$.

\begin{example} $\boldsymbol{\mathcal{Q}}^{\mathrm
h}(2)=\conv\bigl(\pmb{0},(1,-2,1), (0,1,0),(1,-1,1)\bigr)$, see
Figure~{\rm\ref{figure:1}(b)}.
\end{example}

Another description of polytope~(\ref{eq:10}) is
\begin{equation*}
\boldsymbol{\mathcal{Q}}^{\mathrm{h}}(m)=
\lin\bigl(\pmb{h}(\overline{\mathbf{2}^{[1]}};m),\
\pmb{h}(\overline{\mathbf{2}^{[3]}};m),\ \ldots,\
\pmb{h}(\overline{\mathbf{2}^{[m-1]}};m),\ \boldsymbol{\iota}(m)
\bigr)\ \cap\ \mathrm{S}_m\bigl(\boldsymbol{\varPi}(m)\bigr)
\end{equation*}
in the case of $m$ even, and
\begin{equation*}
\boldsymbol{\mathcal{Q}}^{\mathrm{h}}(m)=
\lin\bigl(\pmb{h}(\overline{\mathbf{2}^{[2]}};m),\
\pmb{h}(\overline{\mathbf{2}^{[4]}};m),\ \ldots,\
\pmb{h}(\overline{\mathbf{2}^{[m-1]}};m),\ \boldsymbol{\iota}(m)
\bigr)\ \cap\ \mathrm{S}_m\bigl(\boldsymbol{\varPi}(m)\bigr)\ ,
\end{equation*}
if $m$ is odd. In other words, for any $m$, we have
\begin{equation*}
\boldsymbol{\mathcal{Q}}^{\mathrm{h}}(m) =
\boldsymbol{\mathcal{E}}^{\mathrm{h}}(m)\ \cap\
\mathrm{S}_m\bigl(\boldsymbol{\varPi}(m)\bigr)\ .
\end{equation*}

\subsection{Biorthogonality}

By the {\em principle of
biorthogonality}~\cite[Theorem~1.4.7]{HJ}, we have
\begin{equation*}
\pmb{h}(\overline{\mathbf{2}^{[s]}};m)\cdot
\pmb{h}(\overline{\mathbf{2}^{[t]}};m)^{\top}=0\ ,
\end{equation*}
for all positive integers $s$ and $t$ less than or equal to $m$
such that $s\not\equiv t\pmod{2}$. Indeed, one of the vectors
$\pmb{h}(\overline{\mathbf{2}^{[s]}};m)$ and
$\pmb{h}(\overline{\mathbf{2}^{[t]}};m)$ is a left eigenvector of
the backward identity matrix $\mathbf{U}(m)$ corresponding to an
eigenvalue $\lambda\in\{-1,1\}$, while the other vector is a right
eigenvector corresponding to the other eigenvalue
$\mu\in\{-1,1\}$.

\subsection{The norms}

For a positive integer $k$ such that $k\leq m$, we have
{\footnotesize
\begin{equation*}
\|\pmb{h}(\overline{\mathbf{2}^{[k]}};m)\|^2
=2\Bigl(\binom{2(m-k)}{m-k}-(-1)^k\binom{2(m-k)}{m}\Bigr)
\end{equation*}
} and {\footnotesize
\begin{equation*}
\|\pmb{f}(\overline{\mathbf{2}^{[k]}};m)\|^2=\binom{2k}{k}-1\ .
\end{equation*}
}

\subsection{Orthogonal projectors} $\quad$

\begin{itemize}
\item
If $m$ is even, then the matrix of the orthogonal projector (see,
e.g.,~\cite[Chapter~3]{Strang} on projectors) into
$\boldsymbol{\mathcal{H}}(m)$, relative either to the standard
basis $\bS_m$ or to the basis $\bH^{\blacktriangledown}_m$, is
{\scriptsize
\begin{multline*}
\begin{pmatrix}
\pmb{h}(\overline{\mathbf{2}^{[1]}};m)^{\top} &
\pmb{h}(\overline{\mathbf{2}^{[3]}};m)^{\top}& \cdots&
\pmb{h}(\overline{\mathbf{2}^{[m-1]}};m)^{\top}\end{pmatrix}
\\
\cdot \Biggl(
\begin{pmatrix}
\pmb{h}(\overline{\mathbf{2}^{[1]}};m) \\
\pmb{h}(\overline{\mathbf{2}^{[3]}};m)\\ \vdots \\
\pmb{h}(\overline{\mathbf{2}^{[m-1]}};m)\end{pmatrix}\cdot
\begin{pmatrix}
\pmb{h}(\overline{\mathbf{2}^{[1]}};m)^{\top} &
\pmb{h}(\overline{\mathbf{2}^{[3]}};m)^{\top}& \cdots&
\pmb{h}(\overline{\mathbf{2}^{[m-1]}};m)^{\top}\end{pmatrix}
\Biggr)^{-1}\cdot
\begin{pmatrix}
\pmb{h}(\overline{\mathbf{2}^{[1]}};m) \\
\pmb{h}(\overline{\mathbf{2}^{[3]}};m)\\ \vdots \\
\pmb{h}(\overline{\mathbf{2}^{[m-1]}};m)\end{pmatrix}\ .
\end{multline*}
} The matrix of the orthogonal projector (relative to the standard
basis $\bS_m$) into the subspace $\boldsymbol{\mathcal{F}}(m)$ is
{\scriptsize
\begin{multline*}
\begin{pmatrix}
\pmb{f}(\overline{\mathbf{2}^{[1]}};m)^{\top} &
\pmb{f}(\overline{\mathbf{2}^{[3]}};m)^{\top}& \cdots&
\pmb{f}(\overline{\mathbf{2}^{[m-1]}};m)^{\top}\end{pmatrix}
\\
\cdot \Biggl(
\begin{pmatrix}
\pmb{f}(\overline{\mathbf{2}^{[1]}};m) \\
\pmb{f}(\overline{\mathbf{2}^{[3]}};m)\\ \vdots \\
\pmb{f}(\overline{\mathbf{2}^{[m-1]}};m)\end{pmatrix}\cdot
\begin{pmatrix}
\pmb{f}(\overline{\mathbf{2}^{[1]}};m)^{\top} &
\pmb{f}(\overline{\mathbf{2}^{[3]}};m)^{\top}& \cdots&
\pmb{f}(\overline{\mathbf{2}^{[m-1]}};m)^{\top}\end{pmatrix}
\Biggr)^{-1}\cdot
\begin{pmatrix}
\pmb{f}(\overline{\mathbf{2}^{[1]}};m) \\
\pmb{f}(\overline{\mathbf{2}^{[3]}};m)\\ \vdots \\
\pmb{f}(\overline{\mathbf{2}^{[m-1]}};m)\end{pmatrix}\ .
\end{multline*}
} The matrices of the orthogonal projectors in the case of $m$ odd
are built in an analogous way.

\item
Let $k$ be a positive integer such that $k\leq m$. The $(i,j)$th
entry of the matrix of the orthogonal projector into the
one-dimensional subspace
$\lin\bigl(\pmb{h}(\overline{\mathbf{2}^{[k]}};m)\bigr)$ of the
space $\mathbb{R}^{m+1}$, relative to either of the bases $\bS_m$
and $\bH^{\blacktriangledown}_m$, is {\footnotesize
\begin{equation*}
\frac{(-1)^{i+j}(\binom{m-k}{i}-(-1)^k\binom{m-k}{i-k})
(\binom{m-k}{j}-(-1)^k\binom{m-k}{j-k})}
{2(\binom{2(m-k)}{m-k}-(-1)^k\binom{2(m-k)}{m})} \ ;
\end{equation*}
} the $(i,j)$th entry of the matrix of the orthogonal projector
into the one-dimensional subspace
$\lin\bigl(\pmb{f}(\overline{\mathbf{2}^{[k]}};m)\bigr)$ of
$\mathbb{R}^{m+1}$, relative to $\bS_m$, is {\footnotesize
\begin{equation*}
\frac{\binom{k}{i}\binom{k}{j}}{\binom{2k}{k}-1}\ .
\end{equation*}
}
\end{itemize}
}


\end{document}